\documentclass{amsart}
\usepackage{amsfonts}

\newtheorem{theorem}{Theorem}[section]
\newtheorem{lemma}[theorem]{Lemma}
\theoremstyle{definition}
\newtheorem{definition}[theorem]{Definition}

\theoremstyle{remark}
\newtheorem{remark}[theorem]{Remark}
\numberwithin{equation}{section}
\theoremstyle{plain}

\newtheorem{corollary}{Corollary}

\newtheorem{proposition}{Proposition}

\copyrightinfo{2001}{enter name of copyright holder}

\begin{document}
\title[Gradient Ricci Solitons]{On the Classification of Gradient Ricci
Solitons}
\author{Peter Petersen}
\address{520 Portola Plaza\\
Dept of Math UCLA\\
Los Angeles, CA 90095}
\email{petersen@math.ucla.edu}
\urladdr{http://www.math.ucla.edu/\symbol{126}petersen}
\thanks{}
\author{William Wylie}
\email{wylie@math.ucla.edu}
\urladdr{http://www.math.ucla.edu/\symbol{126}wylie}
\date{}

\begin{abstract}
We show that the only shrinking gradient solitons with vanishing Weyl tensor
are quotients of the standard ones $S^{n},$ $S^{n-1}\times \mathbb{R},$ and $%
\mathbb{R}^{n}$. This gives a new proof of the Hamilton-Ivey-Perel'man
classification of 3-dimensional shrinking gradient solitons. We also show
that gradient solitons with constant scalar curvature and suitably decaying
Weyl tensor when noncompact are quotients of $\mathbb{H}^{n},$ $\mathbb{H}%
^{n-1}\times \mathbb{R},$ $\mathbb{R}^{n},$ $S^{n-1}\times \mathbb{R},$ or $%
S^{n}.$
\end{abstract}

\maketitle

\section{Introduction}

A Ricci soliton is a Riemannian metric together with a vector field $\left(
M,g,X\right) $ that satisfies%
\begin{equation*}
\mathrm{Ric}+\frac{1}{2}L_{X}g=\lambda g.
\end{equation*}%
It is called shrinking when $\lambda >0,$ steady when $\lambda =0$, and
expanding when $\lambda <0$. In case $X=\nabla f$ the equation can also be
written as%
\begin{equation*}
\mathrm{Ric}+\mathrm{Hess}f=\lambda g
\end{equation*}%
and the metric is called a gradient Ricci soliton.

In dimension $2$ Hamilton proved that the shrinking gradient Ricci solitons
with bounded curvature are $S^{2},$ $\mathbb{R}P^{2},$ and $\mathbb{R}^{2}$
with constant curvature \cite{HamiltonI}. Ivey proved the first
classification result in dimension $3$ showing that compact shrinking
gradient solitons have constant positive curvature \cite{Ivey}. In the
noncompact case Perelman has shown that the $3$-dimensional shrinking
gradient Ricci solitons with bounded nonnegative sectional curvature are $%
S^{3}$, $S^{2}\times \mathbb{R}$ and $\mathbb{R}^{3}$ or quotients thereof 
\cite{PerelmanII}. (In Perelman's paper he also includes the assumption that
the manifold is $\kappa $-noncollapsed but this assumption is not necessary
to the argument, see for example, \cite{Chow-Lu-Ni}.) The Hamilton-Ivey
estimate shows that all $3$-dimensional shrinking Ricci solitons with
bounded curvature have non-negative sectional curvature (\cite{Chow-Lu-Ni},
Theorem 6.44) So the work of Perelman, Hamilton, and Ivey together give the
following classification in dimension 3. \footnote{%
We thank Ben Chow for pointing this out to us.}

\begin{theorem}
\label{shrinker-3dim} The only three dimensional shrinking gradient Ricci
solitons with bounded curvature are the finite quotients of $\mathbb{R}^3$, $%
S^2 \times \mathbb{R}$, and $S^3$.
\end{theorem}

Recently Ni and Wallach \cite{Ni-Wallach} have given an alternative approach
to proving the classification of $3$-dimensional shrinkers which extends to
higher dimensional manifolds with zero Weyl tensor. (Every $3$-manifold has
zero Weyl tensor.) Their argument also requires non-negative Ricci
curvature. Also see Naber's paper \cite{Naber} for a different argument in
the 3-dimensional case. By using a different set of formulas we remove the
non-negative curvature assumption.

\begin{theorem}
\label{shrinker-ndim} Let $(M^{n},g,f)$ be a complete shrinking gradient
Ricci soliton of dimension $n\geq 3$ such that $\int_{M}\left\vert \mathrm{%
Ric}\right\vert ^{2}e^{-f}d\mathrm{vol}_{g}<\infty $ and $W=0$ then $M$ is a
finite quotient of $\mathbb{R}^{n}$, $S^{n-1}\times \mathbb{R}$, or $S^{n}$.
\end{theorem}

Recall that a result of Morgan \cite{Morgan} implies $e^{-f}d\mathrm{vol}_g$
is a finite measure so as a corollary we obtain a new direct proof of
Theorem \ref{shrinker-3dim} that does not require the Hamilton-Ivey
estimate. When $M$ is compact Theorem \ref{shrinker-ndim} was established
using similar techniques by Eminenti, LaNave, and Mantegazza in \cite%
{Eminenti-LaNave-Mantegazza}.

We note that a shrinking soliton has finite fundamental group \cite{Wylie}
and that Naber has shown that it can be made into a gradient soliton by
adding an appropriate Killing field to $X$ \cite{Naber}. In the compact case
this was proven by Perelman \cite{PerelmanI}.

If we relax the Weyl curvature condition and instead assume that the scalar
curvature is constant we also get a nice general classification.

\begin{theorem}
\label{expander-ndim} Let $(M^{n},g,f)$ be a complete gradient Ricci soliton
with $n\geq 3$, constant scalar curvature, and $W\left( \nabla f,\cdot
,\cdot ,\nabla f\right) =o\left( \left\vert \nabla f\right\vert ^{2}\right)
, $ then $M$ is a flat bundle of rank $0,$ $1,$ or $n$ over an Einstein
manifold.
\end{theorem}

Note that the Weyl curvature condition is vacuous when $n=3$ or $M$ is
compact. Moreover, when $M$ is compact or the soliton is steady it is
already known that it has to be rigid when the scalar curvature is constant 
\cite{Petersen-WylieI}. It is also worth pointing out that the theorem is in
a sense optimal. Namely, rigid solitons with zero, one or $n$-dimensional
Euclidean factors have $W\left( \nabla f,\cdot ,\cdot ,\nabla f\right) =0.$
When $n=3$ the theorem yields the following new result.

\begin{corollary}
The only 3-dimensional expanding gradient Ricci solitons with constant
scalar curvature are quotients of $\mathbb{R}^{3},$ $\mathbb{H}^{2}\times 
\mathbb{R}$, and $\mathbb{H}^{3}$.
\end{corollary}

There are now many examples of non-trivial gradient solitons, but we do not
know of any with constant scalar curvature. Moreover, we have shown that any
gradient soliton which is homogeneous or has constant scalar curvature and
is radially flat (i.e. $\mathrm{sec}(\nabla f,E)=0$) is a product of
Einstein and Euclidean manifolds \cite{Petersen-WylieII, Petersen-WylieI}.
There are a number non-trivial homogeneous expanding Ricci solitons, even in
dimension 3 (see \cite{Baird-Danielo, Lauret, Lott}). Unlike the shrinking
case, these metrics do not support a \emph{gradient} soliton structure.

Our results follow from considering elliptic equations for various curvature
quantities on solitons. While there are well-known Ricci flow versions of a
number of these formulas, the elliptic proofs are surprisingly
straight-forward and give some interesting extra rigidity. For example, by
considering the equation for the curvature operator we show that if the 
\emph{second} eigenvalue of the curvature operator of a shrinking gradient
Ricci soliton is nonnegative then the metric has nonnegative curvature
operator. This then extends a number of rigidity theorems for nonnegative(or
2-nonnegative) curvature operator (see \cite{Bohm-Wilking, Cao, Naber,
Petersen-WylieII}).

The paper is organized as follows. We start by deriving the formulas for $f$%
-Laplacians of various functions and tensors related to curvature. Next we
give our constant scalar curvature characterization. The proof of this
result is almost entirely algebraic. By contrast the proof of the theorems
for shrinking solitons relies far more heavily on analytic techniques. In
the appendix we review the proof of the classification of 2-dimensional
solitons giving a proof that follows from an Obata type characterization of
warped product manifolds found in \cite{Cheeger-Colding} that does not seem
to appear elsewhere in the literature.

\section{The $f$-Laplacian of Curvature}

In this section we are interested in deriving elliptic equations for the
curvature of a gradient Ricci soliton. Let $V$ be a tensor bundle on a
gradient Ricci soliton, $\nabla$ be the Riemannian connection on $V$, and $T$
a self-adjoint operator on $V$. The $X$-Laplacian and $f$-Laplacian of $T$
is the operator 
\begin{eqnarray*}
(\Delta _{X}T) &=&(\Delta T)-(\nabla _{X}T) \\
(\Delta _{f}T) &=&(\Delta T)-(\nabla _{\nabla f}T)
\end{eqnarray*}%
Where $\Delta$ is the connection Laplacian induced by $\nabla$. We are
interested in the cases where $V = \wedge^2 M$ and $T$ is the curvature
operator $\mathcal{R}$ and where $V=TM$ and $T$ is the $(1,1)$ Ricci tensor.

The formulas in these cases are the following.

\begin{lemma}
\label{CurvForm} For a gradient Ricci soliton 
\begin{eqnarray*}
\Delta _{f}\mathcal{R} &=&2\lambda \mathcal{R}-2\left(\mathcal{R}^{2}+\mathcal{R}%
^{\#}\right) \\
\Delta _{f}\mathrm{Ric} &=&2\lambda \mathrm{Ric}-2\sum_{i=1}^{n}R(\cdot
,E_{i})(\mathrm{Ric}(E_{i})) \\
\Delta _{f}\mathrm{scal} &=&2\lambda \mathrm{scal}-2|\mathrm{Ric}|^{2}
\end{eqnarray*}
\end{lemma}
\begin{remark} A mild warning is in order for the first equation.  We define the induced metric on  $\wedge^2 M$  so that if $\{E_i\}$ is an orthonormal basis of $T_pM$ then $\{ E_i \wedge E_j \}_{i<j}$ is an orthonormal basis of $\wedge^2 T_pM$.  This convention agrees with \cite{Bohm-Wilking} but differs from \cite{Chow-Knopf, HamiltonII, HamiltonIII}.
\end{remark}
\begin{remark}
The last equation is well known, see \cite{Petersen-WylieI} for a proof. A
similar equation for the Ricci tensor appears in \cite%
{Eminenti-LaNave-Mantegazza}. Some other interesting formulas for the
curvature operator of gradient solitons appear in \cite{Cao}.
\end{remark}

\begin{remark}
Since Ricci solitons are special solutions to the Ricci flow the above
equations can be derived from the parabolic formulas derived by Hamilton 
\cite{HamiltonII} for the Ricci flow 
\begin{eqnarray*}
\frac{\partial }{\partial t}\mathcal{R} &=&\Delta \mathcal{R}+2 \left(\mathcal{R}%
^{2}+\mathcal{R}^{\#}\right) \\
\frac{\partial }{\partial t}\mathrm{Ric} &=&\Delta _{L}\mathrm{Ric} \\
\frac{\partial }{\partial t}\mathrm{scal} &=&\Delta \mathrm{scal}+2|\mathrm{%
Ric}|^{2}.
\end{eqnarray*}%
However, we will simply perform the elliptic calculation which is more
straight forward (for example no \textquotedblleft Uhlenbeck trick" is
necessary). We also expect similar calculations will give formulas for
elliptic equations which do not come directly from a Ricci flow.
\end{remark}

$\mathcal{R}^{\#}$ is the Lie-Algebra square of $\mathcal{R}$ introduced by
Hamilton in \cite{HamiltonIII}. Recall that if we change $\mathcal{R}^{\#}$
into a (0,4)-tensor its formula is 
\begin{eqnarray*}
g\left( \mathcal{R}^{\#}(X\wedge Y),W\wedge Z\right) &=&R^{\#}(X,Y,Z,W) \\
&=& B(X,W,Y,Z)-B(X,Z,Y,W),
\end{eqnarray*}%
where 
\begin{equation*}
B(X,Y,W,Z)=-\sum_{i=1}^{n}g\left( R(X,E_{i})Y,R(W,E_{i})Z\right)
\end{equation*}
and $\{E_{i}\}$ is an orthonrmal basis of $T_{p}M.$ It is also convenient to
identify $\wedge ^{2}T_{p}M$ with $\mathfrak{so}(n)$. Then $\wedge
^{2}T_{p}M $ becomes a Lie Algebra and the formula for $\mathcal{R}^{\#}$
becomes 
\begin{equation*}
g(\mathcal{R}^{\#}(U),V)=\frac{1}{2} \sum_{\alpha ,\beta }g\left( [\mathcal{R}(\phi
_{\alpha }),\mathcal{R}(\phi _{\beta })],U\right) g\left( [\phi _{\alpha
},\phi _{\beta }],V\right)
\end{equation*}%
for any two bi-vectors $U$ and $V$, where $\{\phi _{\alpha }\}$ is an
orthnormal basis of $\wedge ^{2}T_{p}M$. (see page 186 of \cite{Chow-Knopf}
for the derivation of the equivalence of these two formulas. See \cite%
{Bohm-Wilking, Chen} for more about $\mathcal{R}^{\#}$.)

Before the main calculation we recall some curvature identities for gradient
Ricci solitons.

\begin{proposition}
For a gradient Ricci soliton 
\begin{equation}
\nabla \mathrm{scal}=2\mathrm{div}\left( \mathrm{Ric}\right) =2\mathrm{Ric}%
\left( \nabla f\right)
\end{equation}%
\begin{equation}
(\nabla _{X}\mathrm{Ric})(Y)-(\nabla _{Y}\mathrm{Ric})(X)=-R(X,Y)\nabla f
\label{two}
\end{equation}%
\begin{equation}
\nabla _{\nabla f}\mathrm{Ric}+\mathrm{Ric}\circ \left( \lambda I-\mathrm{Ric%
}\right) =R\left( \cdot ,\nabla f\right) \nabla f+\frac{1}{2}\nabla _{\cdot
}\nabla \mathrm{scal}  \label{three}
\end{equation}%
\begin{equation}
\sum_{i=1}^{n}(\nabla _{E_{i}}R)(E_{i},X,Y)=R(\nabla f,X)Y  \label{four}
\end{equation}
\end{proposition}

\begin{proof}
The proofs of the first three identities can be found in \cite%
{Petersen-WylieI}. For the fourth formula consider 
\begin{eqnarray*}
\sum_{i=1}^{n}(\nabla _{E_{i}}R)(E_{i},X,Y,Z) &=&-\sum_{i=1}^{n}(\nabla
_{E_{i}}R)(Y,Z,X,E_{i}) \\
&=&-(\mathrm{div}R)(Y,Z,X) \\
&=&-\left( \nabla _{Y}\mathrm{Ric}\right) \left( Z,X\right) +\left( \nabla
_{Z}\mathrm{Ric}\right) \left( Y,X\right) \\
&=&g(R(Y,Z)\nabla f,X) \\
&=&g(R(\nabla f,X)Y,Z)
\end{eqnarray*}%
Where in the third line we have used the (contracted) 2nd Bianchi identity
and in the fourth line we have used (\ref{two}).
\end{proof}

We are now ready to derive the formula for the $f$-Laplacian of curvature.

\begin{proof}[Proof of Lemma \protect\ref{CurvForm}]
We will begin the calculation by considering the (0,4)-curvature tensor $R$.
Fix a point $p$, let $X,Y,Z,W$ be vector fields with $\nabla X=\nabla
Y=\nabla Z=\nabla W=0$ at $p$ and let $E_{i}$ be normal coordinates at $p$.
Then 
\begin{eqnarray*}
\left( \Delta R\right) \left( X,Y,Z,W\right) &=&\sum_{i=1}^{n}\left(
\nabla^2 _{E_{i},E_{i}}R\right) \left( X,Y,Z,W\right) \\
&=&\sum_{i=1}^{n}\left( \nabla^2 _{E_{i},X}R\right) \left(
E_{i},Y,Z,W\right) -\left( \nabla^2 _{E_{i},Y}R\right) \left(
E_{i},X,Z,W\right) \\
&=&\sum_{i=1}^{n}\left( \nabla^2 _{X,E_{i}}R\right) \left(
E_{i},Y,Z,W\right) -\left( \nabla^2 _{Y,E_{i}}R\right) \left(
E_{i},X,Z,W\right) \\
&&\qquad +\left( R_{E_{i},X}R\right) \left( E_{i},Y,Z,W\right) -\left(
R_{E_{i},Y}R\right) \left( E_{i},X,Z,W\right) \\
&=&\nabla _{X}\left( R\left( \nabla f,Y,Z,W\right) \right) -\nabla
_{Y}\left( R\left( \nabla f,X,Z,W\right) \right) \\
&&+\sum_{i=1}^{n}\left( R_{E_{i},X}R\right) \left( E_{i},Y,Z,W\right)
-\left( R_{E_{i},Y}R\right) \left( E_{i},X,Z,W\right) \\
&=&\left( \nabla _{X}R\right) \left( \nabla f,Y,Z,W\right) +R\left( \nabla
_{X}\nabla f,Y,Z,W\right) \\
&&-\left( \nabla _{Y}R\right) \left( \nabla f,X,Z,W\right) -R\left( \nabla
_{Y}\nabla f,X,Z,W\right) \\
&&+\sum_{i=1}^{n}\left( R_{E_{i},X}R\right) \left( E_{i},Y,Z,W\right)
-\left( R_{E_{i},Y}R\right) \left( E_{i},X,Z,W\right)
\end{eqnarray*}

Where in the fourth line we have applied (\ref{four}). The second Bianchi
identity implies 
\begin{equation*}
\left( \nabla _{\nabla f}R\right) \left( X,Y,Z,W\right) =\left( \nabla
_{X}R\right) \left( \nabla f,Y,Z\right) -\left( \nabla _{Y}R\right) \left(
\nabla f,X,Z\right) ,
\end{equation*}%
and the gradient soliton equation gives 
\begin{equation*}
R\left( \nabla _{X}\nabla f,Y,Z,W\right) =\lambda R\left( X,Y,Z,W\right)
-R\left( \mathrm{Ric}(X),Y,Z,W\right) .
\end{equation*}%
So we have 
\begin{eqnarray*}
\left( \Delta _{f}R\right) (X,Y,Z,W) &=&2\lambda R\left( X,Y,Z,W\right)
-R\left( \mathrm{Ric}(X),Y,Z,W\right) +R\left( \mathrm{Ric}(Y),X,Z,W\right) 
\\
&&+\sum_{i=1}^{n}\left( R_{E_{i},X}R\right) \left( E_{i},Y,Z,W\right)
-\left( R_{E_{i},Y}R\right) \left( E_{i},X,Z,W\right) 
\end{eqnarray*}%
We now must unravel the the terms remaining inside the sum. By definition 
\begin{eqnarray*}
\left( R_{E_{i},X}R\right) \left( E_{i},Y,Z,W\right)  &=&R\left(
E_{i},X,R(E_{i},Y)Z,W\right) -R\left( R(E_{i},X)E_{i},Y,Z,W\right)  \\
&&-R\left( E_{i},R(E_{i},X)Y,Z,W\right) -R\left(
E_{i},Y,R(E_{i},X)Z,W\right) .
\end{eqnarray*}%
A straight forward calculation involving the Bianchi identity then gives 
\begin{eqnarray*}
&&\sum_{i=1}^{n}\left( R_{E_{i},X}R\right) \left( E_{i},Y,Z,W\right) -\left(
R_{E_{i},Y}R\right) \left( E_{i},X,Z,W\right)  \\
&=&R\left( \mathrm{Ric}(X),Y,Z,W\right) -R\left( \mathrm{Ric}%
(Y),X,Z,W\right)  \\
&&+\sum_{i=1}^{n}\left[ -2R\left( X,E_{i},R(E_{i},Y)Z,W\right) +2R\left(
Y,E_{i},R(E_{i},X)Z,W\right) \right]  \\
&&+\sum_{i=1}^{n}R\left( E_{i},R(X,Y)E_{i},Z,W\right) 
\end{eqnarray*}%
We now have 
\begin{eqnarray*}
&&\left( \Delta _{f}R\right) (X,Y,Z,W) \\
&=&2\lambda R\left( X,Y,Z,W\right) +2\sum_{i=1}^{n}\left[ -R\left(
X,E_{i},R(E_{i},Y)Z,W\right) +R\left( Y,E_{i},R(E_{i},X)Z,W\right) \right] 
\\
&&+\sum_{i=1}^{n}R\left( E_{i},R(X,Y)E_{i},Z,W\right) .
\end{eqnarray*}%
However, 
\begin{eqnarray*}
\sum_{i=1}^{n}R\left( E_{i},R(X,Y)E_{i},Z,W\right) 
&=&-\sum_{i=1}^{n}R\left( W,Z,E_{i},R(X,Y)E_{i}\right)  \\
&=&-\sum_{i=1}^{n}g\left( R(W,Z)E_{i},R(X,Y)E_{i}\right)  \\
&=&-2R^{2}(X,Y,Z,W)
\end{eqnarray*}%
and 
\begin{eqnarray*}
&&2\sum_{i=1}^{n}-R\left( X,E_{i},R(E_{i},Y)Z,W\right) +R\left(
Y,E_{i},R(E_{i},X)Z,W\right)  \\
&=&2\sum_{i=1}^{n}-g\left( R(X,E_{i})W,R(Y,E_{i})Z\right) +g\left(
R(Y,E_{i})W,R(X,E_{i})Z\right)  \\
&=&-2R^{\#}(X,Y,Z,W).
\end{eqnarray*}%
So we have obtained the desired formula for the curvature operator.

To compute the formula for the Ricci tensor we could trace the formula for
the curvature tensor, or we can give the following direct proof.

Again fix a point $p$, extend $Y(p)$ to a vector field in a neighborhood of $%
p$ such that $\nabla Y=0$, and let $E_{i}$ be normal coordinates at $p$,
then 
\begin{eqnarray*}
(\Delta \mathrm{Ric})(Y) &=&\sum_{i=1}^{n}\left( \nabla _{E_{i},E_{i}}^{2}%
\mathrm{Ric}\right) (Y) \\
&=&\sum_{i=1}^{n}\left( (\nabla _{E_{i},Y}^{2}\mathrm{Ric})(E_{i})-\nabla
_{E_{i}}\left( R(E_{i},Y)\nabla f\right) \right) \\
&=&\sum_{i=1}^{n}\left( (\nabla _{Y,E_{i}}^{2}\mathrm{Ric}%
)(E_{i})-(R_{Y,E_{i}}\mathrm{Ric})(E_{i})-(\nabla _{E_{i}}R)(E_{i},Y,\nabla
f)-R(E_{i},Y)(\nabla _{E_{i}}\nabla f)\right) \\
&=&\nabla _{Y}(\mathrm{div}(\mathrm{Ric}))+\mathrm{Ric}(\mathrm{Ric}%
(Y))+R(Y,\nabla f)\nabla f+\lambda \mathrm{Ric}(Y)-2\sum_{i=1}^{n}R(Y,E_{i})(%
\mathrm{Ric}(E_{i})) \\
&=&\left( \nabla _{\nabla f}\mathrm{Ric}\right) \left( Y\right) +2\lambda 
\mathrm{Ric}(Y)-2\sum_{i=1}^{n}R(Y,E_{i})(\mathrm{Ric}(E_{i})).
\end{eqnarray*}%
Where in going from the first to second lines we have applied (\ref{two}),
in going from the third to fourth lines we apply (\ref{four}) and in
obtaining the last line we apply (\ref{three}).
\end{proof}

From the Ricci equation we can also derive the following formula

\begin{lemma}
\begin{eqnarray*}
\Delta _{f}\left( \mathrm{Ric}(\nabla f,\nabla f)\right) &=&4\lambda \mathrm{%
Ric}(\nabla f,\nabla f)-2D_{\nabla f}\left\vert \mathrm{Ric}\right\vert ^{2}
\\
&&+2\mathrm{Ric}\left( \nabla _{E_{i}}\nabla f,\nabla _{E_{i}}\nabla
f\right) +2\sum_{i=1}^{n}R(\nabla f,E_{i},\mathrm{Ric}(E_{i}),\nabla f)
\end{eqnarray*}%
or equivalently%
\begin{equation*}
\frac{1}{2}\Delta _{f}\left( D_{\nabla f}\mathrm{scal}\right) =D_{\nabla
f}\Delta _{f}\mathrm{scal}+2\mathrm{Ric}\left( \nabla _{E_{i}}\nabla
f,\nabla _{E_{i}}\nabla f\right) +2\sum_{i=1}^{n}R(\nabla f,E_{i},\mathrm{Ric%
}(E_{i}),\nabla f)
\end{equation*}
\end{lemma}

\begin{proof}
From the above equation we get%
\begin{eqnarray*}
\Delta _{f}\left( \mathrm{Ric}(\nabla f,\nabla f)\right) &=&\left( \Delta
_{f}\mathrm{Ric}\right) \left( \nabla f,\nabla f\right) +2\mathrm{Ric}%
(\Delta _{f}\nabla f,\nabla f)+2\mathrm{Ric}\left( \nabla _{E_{i}}\nabla
f,\nabla _{E_{i}}\nabla f\right) \\
&&+4\left( \nabla _{E_{i}}\mathrm{Ric}\right) \left( \nabla _{E_{i}}\nabla
f,\nabla f\right) \\
&=&\left( \Delta _{f}\mathrm{Ric}\right) \left( \nabla f,\nabla f\right)
-2\lambda \mathrm{Ric}(\nabla f,\nabla f)+2\mathrm{Ric}\left( \nabla
_{E_{i}}\nabla f,\nabla _{E_{i}}\nabla f\right) \\
&&+4\lambda \left( \nabla _{E_{i}}\mathrm{Ric}\right) \left( E_{i},\nabla
f\right) -4\left( \nabla _{E_{i}}\mathrm{Ric}\right) \left( \mathrm{Ric}%
\left( E_{i}\right) ,\nabla f\right) \\
&=&-2R(\nabla f,E_{i},\mathrm{Ric}(E_{i}),\nabla f)+2\mathrm{Ric}\left(
\nabla _{E_{i}}\nabla f,\nabla _{E_{i}}\nabla f\right) \\
&&+4\lambda \mathrm{Ric}\left( \nabla f,\nabla f\right) -4\left( \nabla
_{\nabla f}\mathrm{Ric}\right) \left( \mathrm{Ric}\left( E_{i}\right)
,E_{i}\right) +4R(\nabla f,E_{i},\mathrm{Ric}(E_{i}),\nabla f) \\
&=&4\lambda \mathrm{Ric}\left( \nabla f,\nabla f\right) -2D_{\nabla
f}\left\vert \mathrm{Ric}\right\vert ^{2} \\
&&+2\mathrm{Ric}\left( \nabla _{E_{i}}\nabla f,\nabla _{E_{i}}\nabla
f\right) +2R(\nabla f,E_{i},\mathrm{Ric}(E_{i}),\nabla f)
\end{eqnarray*}

The second formula follows from 
\begin{eqnarray*}
\Delta _{f}\mathrm{scal} &=&2\lambda \mathrm{scal}-2\left\vert \mathrm{Ric}%
\right\vert ^{2}, \\
2\mathrm{Ric}\left( \nabla f,\nabla f\right) &=&D_{\nabla f}\mathrm{scal}
\end{eqnarray*}
\end{proof}

We are now going to see how the Weyl decomposition affects the formula for
the Ricci tensor.

\begin{lemma}
\begin{eqnarray*}
\Delta _{f}\mathrm{Ric} &=&2\lambda \mathrm{Ric}-\frac{2n\mathrm{scal}}{%
\left( n-1\right) \left( n-2\right) }\mathrm{Ric}+\frac{4}{n-2}\mathrm{Ric}%
^{2} \\
&&-\frac{2}{\left( n-2\right) }\left( \left\vert \mathrm{Ric}\right\vert
^{2}-\frac{\mathrm{scal}^{2}}{n-1}\right) I+W\left( \cdot ,E_{i},\mathrm{Ric}%
\left( E_{i}\right) \right)
\end{eqnarray*}%
and%
\begin{eqnarray*}
\frac{1}{2}\Delta _{f}\left( D_{\nabla f}\mathrm{scal}\right) &=&D_{\nabla
f}\Delta _{f}\mathrm{scal}+2\mathrm{Ric}\left( \nabla _{E_{i}}\nabla
f,\nabla _{E_{i}}\nabla f\right) \\
&&+\frac{2n\mathrm{scal}}{\left( n-1\right) \left( n-2\right) }\mathrm{Ric}%
\left( \nabla f,\nabla f\right) -\frac{4}{n-2}\mathrm{Ric}\left( \mathrm{Ric}%
\left( \nabla f\right) ,\nabla f\right) \\
&&+\frac{2}{\left( n-2\right) }\left( \left\vert \mathrm{Ric}\right\vert
^{2}-\frac{\mathrm{scal}^{2}}{n-1}\right) \left\vert \nabla f\right\vert
^{2}+W\left( \nabla f,E_{i},\mathrm{Ric}\left( E_{i}\right) ,\nabla f\right)
\end{eqnarray*}
\end{lemma}

\begin{proof}
The Weyl decomposition looks like%
\begin{eqnarray*}
R &=&W+\frac{1}{n-2}\mathrm{Ric}\circ g-\frac{\mathrm{scal}}{2\left(
n-1\right) \left( n-2\right) }g\circ g, \\
h\circ g\left( x,y,y,x\right) &=&h\left( x,x\right) g\left( y,y\right)
+h\left( y,y\right) g\left( x,x\right) -2h\left( x,y\right) g\left(
x,y\right)
\end{eqnarray*}%
where $W$ is absent when $n=3.$ More specifically we need%
\begin{eqnarray*}
R\left( x,y,y,x\right) &=&\frac{1}{n-2}\left( \mathrm{Ric}\left( x,x\right)
g\left( y,y\right) +\mathrm{Ric}\left( y,y\right) g\left( x,x\right) -2%
\mathrm{Ric}\left( x,y\right) g\left( x,y\right) \right) \\
&&-\frac{\mathrm{scal}}{\left( n-1\right) \left( n-2\right) }\left(
\left\vert x\right\vert ^{2}\left\vert y^{2}\right\vert -\left( g\left(
x,y\right) \right) ^{2}\right) +W\left( x,y,y,x\right)
\end{eqnarray*}%
If we assume that $E_{i}$ is an orthonormal frame that diagonalizes the
Ricci tensor $\mathrm{Ric}\left( E_{i}\right) =\rho _{i}E_{i},$ then

\begin{eqnarray*}
\mathrm{Ric}\left( Y\right) &=&g\left( Y,E_{i}\right) \rho _{i}E_{i} \\
\mathrm{Ric}\left( Y,Y\right) &=&\rho _{i}\left( g\left( Y,E_{i}\right)
\right) ^{2} \\
\mathrm{Ric}\left( Y,\mathrm{Ric}\left( Y\right) \right) &=&\mathrm{Ric}%
\left( Y,g\left( Y,E_{i}\right) \rho _{i}E_{i}\right)
\end{eqnarray*}%
using this the Weyl free part of the formula for 
\begin{equation*}
R(Y,E_{i},\mathrm{Ric}(E_{i}),Y)=\rho _{i}R(Y,E_{i},E_{i},Y)
\end{equation*}%
becomes%
\begin{eqnarray*}
&&\frac{1}{n-2}\left( \mathrm{scal}\cdot \mathrm{Ric}\left( Y,Y\right) +\rho
_{i}^{2}\left\vert Y\right\vert ^{2}-2\mathrm{Ric}\left( Y,\rho
_{i}E_{i}\right) g\left( Y,E_{i}\right) \right) \\
&&-\frac{\mathrm{scal}}{\left( n-1\right) \left( n-2\right) }\left(
\left\vert Y\right\vert ^{2}\mathrm{scal}-\rho _{i}\left( g\left(
Y,E_{i}\right) \right) ^{2}\right) \\
&=&\frac{1}{n-2}\left( \mathrm{scal}\cdot \mathrm{Ric}\left( Y,Y\right)
+\left\vert \mathrm{Ric}\right\vert ^{2}\left\vert Y\right\vert ^{2}-2%
\mathrm{Ric}\left( Y,\mathrm{Ric}\left( Y\right) \right) \right) \\
&&+\frac{\mathrm{scal}}{\left( n-1\right) \left( n-2\right) }\left( \mathrm{%
Ric}\left( Y,Y\right) -\left\vert Y\right\vert ^{2}\mathrm{scal}\right) \\
&=&\frac{1}{n-2}\left( \left\vert \mathrm{Ric}\right\vert ^{2}-\frac{\mathrm{%
scal}^{2}}{n-1}\right) \left\vert Y\right\vert ^{2}-\frac{2}{n-2}\mathrm{Ric}%
\left( Y,\mathrm{Ric}\left( Y\right) \right) \\
&&+\left( \frac{\mathrm{scal}}{\left( n-1\right) \left( n-2\right) }+\frac{%
\mathrm{scal}}{\left( n-2\right) }\right) \mathrm{Ric}\left( Y,Y\right) \\
&=&\frac{1}{n-2}\left( \left\vert \mathrm{Ric}\right\vert ^{2}-\frac{\mathrm{%
scal}^{2}}{n-1}\right) \left\vert Y\right\vert ^{2}-\frac{2}{n-2}\mathrm{Ric}%
\left( Y,\mathrm{Ric}\left( Y\right) \right) \\
&&+\frac{n\mathrm{scal}}{\left( n-1\right) \left( n-2\right) }\mathrm{Ric}%
\left( Y,Y\right)
\end{eqnarray*}%
This establishes the first formula and the second by using $Y=\nabla f.$
\end{proof}

\section{Constant Scalar Curvature}

We now turn our attention to the case where $\mathrm{scal}$ is constant in
dimensions $n\geq 3$.

Recall the following results from \cite{Petersen-WylieI} (Propositions 5 and
7).

\begin{proposition}
\label{Prop} Assume that we have a shrinking (resp. expanding) gradient
soliton%
\begin{equation*}
\mathrm{Ric}+\mathrm{Hess}f=\lambda g
\end{equation*}%
with constant scalar curvature. Then $0\leq \mathrm{scal}\leq n\lambda $
(resp. $n\lambda \leq \mathrm{scal}\leq 0$.) Moreover, the metric is flat
when $\mathrm{scal}=0$ and Einstein when $\mathrm{scal}=n\lambda $. In
addition $f$ is unbounded when $M$ is noncompact and $\mathrm{scal}\neq
n\lambda $.
\end{proposition}

This in conjunction with the above formulas allow us to prove

\begin{theorem}
Any gradient soliton with constant scalar curvature, $\lambda \neq 0$ and $%
W\left( \nabla f,\cdot ,\cdot ,\nabla f\right) =o\left( \left\vert \nabla
f\right\vert ^{2}\right) $ is rigid.
\end{theorem}

\begin{proof}
We can assume that $M$ is noncompact and that $f$ is unbounded. The fact
that the scalar curvature is constant in addition shows that%
\begin{equation*}
0=\Delta _{f}\mathrm{scal}=\lambda \mathrm{scal}-\left\vert \mathrm{Ric}%
\right\vert ^{2}
\end{equation*}%
and from the formula for $\Delta _{f}\mathrm{Ric}\left( \nabla f,\nabla
f\right) $ we get 
\begin{eqnarray*}
0 &=&\frac{1}{2}\Delta _{f}\left( D_{\nabla f}\mathrm{scal}\right) \\
&=&D_{\nabla f}\Delta _{f}\mathrm{scal}+2\mathrm{Ric}\left( \nabla
_{E_{i}}\nabla f,\nabla _{E_{i}}\nabla f\right) +2W\left( \nabla f,E_{i},%
\mathrm{Ric}\left( E_{i}\right) ,\nabla f\right) \\
&&+\frac{2n\mathrm{scal}}{\left( n-1\right) \left( n-2\right) }\mathrm{Ric}%
\left( \nabla f,\nabla f\right) -\frac{4}{n-2}\mathrm{Ric}\left( \mathrm{Ric}%
\left( \nabla f\right) ,\nabla f\right) \\
&&+\frac{2}{\left( n-2\right) }\left( \left\vert \mathrm{Ric}\right\vert
^{2}-\frac{\mathrm{scal}^{2}}{n-1}\right) \left\vert \nabla f\right\vert ^{2}
\\
&=&2\mathrm{Ric}\left( \nabla _{E_{i}}\nabla f,\nabla _{E_{i}}\nabla
f\right) +2W\left( \nabla f,E_{i},\mathrm{Ric}\left( E_{i}\right) ,\nabla
f\right) +\frac{2}{\left( n-2\right) }\left( \left\vert \mathrm{Ric}%
\right\vert ^{2}-\frac{\mathrm{scal}^{2}}{n-1}\right) \left\vert \nabla
f\right\vert ^{2}
\end{eqnarray*}%
Since $\left\vert \mathrm{Ric}\right\vert ^{2}=\lambda \mathrm{scal}$ is
constant we see that both $\mathrm{Ric}$ and $\mathrm{Hess}f$ are bounded.
Thus $\mathrm{Ric}\left( \nabla _{E_{i}}\nabla f,\nabla _{E_{i}}\nabla
f\right) $ is bounded and $W\left( \nabla f,E_{i},\mathrm{Ric}\left(
E_{i}\right) ,\nabla f\right) =o\left( \left\vert \nabla f\right\vert
^{2}\right) .$ Recall that%
\begin{equation*}
\mathrm{scal}+\left\vert \nabla f\right\vert ^{2}-2\lambda f=\mathrm{const}
\end{equation*}%
so if the scalar curvature is constant and $f$ is unbounded we see that $%
\left\vert \nabla f\right\vert ^{2}$ is unbounded. This implies that%
\begin{equation*}
\left\vert \mathrm{Ric}\right\vert ^{2}-\frac{\mathrm{scal}^{2}}{n-1}=0
\end{equation*}%
as it is constant. We know in addition that $\mathrm{Ric}$ has one zero
eigenvalue when $\nabla f\neq 0,$ so in that case the Cauchy-Schwarz
inequality shows that 
\begin{equation*}
\left\vert \mathrm{Ric}\right\vert ^{2}\geq \frac{\mathrm{scal}^{2}}{n-1}
\end{equation*}%
with equality holding only if all the other eigenvalues are the same.

If $\nabla f$ vanishes on an open set, then the metric is Einstein on that
set, in particular $\mathrm{scal}=n\lambda $ everywhere and so the entire
metric is Einstein. This means that we can assume $\nabla f\neq 0$ on an
open dense set. Thus $\mathrm{Ric}$ has a zero eigenvalue everywhere and the
other eigenvalues are given by the constant%
\begin{equation*}
\rho =\frac{\mathrm{scal}}{n-1}.
\end{equation*}%
But by Corollary \ref{splittingI} which we will prove below this implies
that $\tilde{M}=N^{n-1}\times \mathbb{R}$ where $N$ is Einstein if $n>3$.
When $n=3$, $N$ is a surface and so must also have constant curvature if $M$
does.
\end{proof}

We now prove that a gradient Ricci soliton whose Ricci curvature has one
nonzero eigenvalue of multiplicity $n-1$ at every point must split. This
will follow from the following more general lemma.

\begin{lemma}
\label{T-lemma} Let $T$ be a constant rank, symmetric, nonnegative tensor on
some (tensor) bundle. If $g\left( \left( \Delta _{X}T\right) \left( s\right)
,s\right) \leq 0$ for $s\in \mathrm{ker}T,$ then the kernel is a parallel
subbundle.
\end{lemma}

\begin{proof}
We are assuming that $\mathrm{ker}T$ is a subbundle. Select an orthonormal
frame $E_{1},...,E_{n}$ and let $s$ be section of $\mathrm{ker}T.$ First
note that%
\begin{equation*}
(\Delta _{X}T)(s)=\Delta _{X}(T(s))-2\sum_{i=1}^{n}\left( \left( \nabla
_{E_{i}}T\right) (\nabla _{E_{i}}s)\right) +T(\Delta _{X}s)
\end{equation*}%
so from the hypothesis we have 
\begin{eqnarray*}
0 &\geq &g((\Delta _{X}T)(s),s) \\
&=&-2\sum_{i=1}^{n}g(\left( \nabla _{E_{i}}T\right) (\nabla
_{E_{i}}s),s)+g(T(\Delta _{X}s),s) \\
&=&-2\sum_{i=1}^{n}g(\nabla _{E_{i}}s,\left( \nabla _{E_{i}}T\right)
(s))+g(\Delta _{X}s,T(s)) \\
&=&-2\sum_{i=1}^{n}g(\nabla _{E_{i}}s,\left( \nabla _{E_{i}}T\right) (s)) \\
&=&2\sum_{i=1}^{n}g(\nabla _{E_{i}}s,T(\nabla _{E_{i}}s))
\end{eqnarray*}%
The nonnegativity of $T$ then gives that $\nabla s\in \mathrm{ker}T$.
\end{proof}

\begin{corollary}
\label{splittingI} Let $(M,g,f)$ be a gradient Ricci soliton such that, at
each point, the Ricci tensor has one nonzero eigenvalue of multiplicity $n-1$%
, then $\tilde{M}=N^{n-1}\times \mathbb{R}.$ Moreover, if $n>3$ then $N$ is
Einstein.
\end{corollary}

\begin{proof}
Let $E_{1},...,E_{n}$ be an orthonormal frame such that $\mathrm{Ric}\left(
E_{1}\right) =0$ and $\mathrm{Ric}\left( E_{i}\right) =\rho E_{i}$ for $i>1.$
Then%
\begin{eqnarray*}
\left( \Delta _{f}\mathrm{Ric}\right) \left( Y\right) &=&2\lambda \mathrm{Ric%
}\left( Y\right) -2\sum_{i=1}^{n}R\left( Y,E_{i}\right) \mathrm{Ric}\left(
E_{i}\right) \\
&=&2\lambda \mathrm{Ric}\left( Y\right) -2\rho \sum_{i=2}^{n}R\left(
Y,E_{i}\right) E_{i} \\
&=&2\left( \lambda -\rho \right) \mathrm{Ric}\left( Y\right) +2\rho R\left(
Y,E_{1}\right) E_{1}.
\end{eqnarray*}%
Since this vanishes on $E_{1}$ we see that the previous lemma can be applied.
\end{proof}

\section{Shrinkers}

For gradient shrinking solitons we use an approach due to Naber (\cite{Naber}%
, section 7). There is a natural measure $e^{-f} d{\mathrm{vol}_g}$ which
makes the $f$-Laplacian self-adjoint. From the perspective of comparison
geometry the tensor $\mathrm{Ric} + \mathrm{Hess}f$ is the Ricci tensor for
this measure and Laplacian. (see e.g. \cite{Lichnerowicz, Morgan, Wei-Wylie}%
). In particular for a shrinking soliton the measure must be bounded above
by a Gaussian measure, note that no assumption on the boundedness of Ricci
curvature is necessary.

\begin{lemma}
\label{vol-comp} (\cite{Morgan}, \cite{Wei-Wylie}) On a shrinking gradient
Ricci soliton the measure $e^{-f}d\mathrm{vol}_{g}$ is finite and if $%
u=O\left( e^{\alpha d^{2}(\cdot ,p)}\right) $ for some $\alpha <\frac{%
\lambda }{2}$ and fixed point $p$ then $u\in L^{2}(e^{-f}d\mathrm{vol}_{g})$.
\end{lemma}

In \cite{Naber} Naber combines a similar volume comparison with a refinement
of a Liouville theorem of Yau (\cite{Yau}, Theorem 3). We will apply the
Liouville theorem to non-smooth functions such as the smallest eigenvalue of
the Ricci tensor so we need to refine these arguments further.

\begin{theorem}[Yau-Naber Liouville Theorem]
Let $\left( M,g,f\right) $ be a manifold with finite $f$-volume: $\int
e^{-f}d\mathrm{vol}<\infty .$ If $u$ is a locally Lipschitz function in $%
L^2(e^{-f} d \mathrm{vol}_g)$ which is bounded below such that 
\begin{equation*}
\Delta_f(u) \geq 0
\end{equation*}
in the sense of barriers, then $u$ is constant.
\end{theorem}

\begin{proof}
Note that since the measure is finite and $u$ is bounded from below we can
assume $u$ is positive by adding a suitable constant to $u$.

To prove the theorem we must modify slightly the techniques of Yau and
Naber. First we apply a heat kernel smoothing procedure of Greene and Wu
(see \cite{Greene-Wu}, section 3).

Let $K$ be a smooth compact subset of $M$ and let $U(x,t)$ be the solution
to the equation 
\begin{eqnarray*}
\left( \frac{\partial }{\partial t}-\Delta _{f}\right) U &=&0 \\
U(x,0) &=&\tilde{u}(x)
\end{eqnarray*}%
on the double of a smooth open set that contains $K$ where $\tilde{u}$ is a
continuous extension of $u$ to the larger open set. Then, by the standard
theory, $U_{t}$ is a smooth function that converges in $W^{1,2}(K)$ to $u$
as $t\rightarrow 0$. Moreover, Green and Wu show that given $\varepsilon >0$
there is $t_{0}$ such that for all $t<t_{0}$ 
\begin{eqnarray*}
\Delta _{f}\left( U(\cdot ,t)\right) &\geq &-\varepsilon ,
\end{eqnarray*}%
on $K$.

Now to the proof of the theorem. Let $x\in M$ and $r_{k}\rightarrow \infty $%
. Using the procedure described above we construct smooth functions $u_{k}$
such that 
\begin{eqnarray*}
|u_{k}-u|_{(W^{1,2}(B(x,r_{k}+1))} &<&\frac{1}{k} \\
\Delta _{f}(u_{k}) &\geq &-\frac{1}{k}
\end{eqnarray*}%
Let $\phi _{k}$ be a cut-off function which is $1$ on $B(x,1)$, $0$ outside
of $B(x,r_{k}+1)$, and has $|\nabla \phi _{k}|\leq \frac{2}{r_{k}}$. First
we integrate by parts.

\begin{equation*}
\int_{M}\Delta _{f}(u_{k})\phi _{k}^{2}u_{k}\left( e^{-f}d\mathrm{vol}%
_{g}\right) =-\int_{M}2\phi _{k}u_{k}g(\nabla u_{k},\nabla \phi _{k})\left(
e^{-f}d\mathrm{vol}_{g}\right) -\int_{M}\phi _{k}^{2}|\nabla
u_{k}|^{2}\left( e^{-f}d\mathrm{vol}_{g}\right) .
\end{equation*}%
Then we complete the square 
\begin{equation*}
\left\vert \sqrt{\frac{1}{2}}\phi _{k}\nabla u_{k}+\sqrt{2}u_{k}\nabla \phi
_{k}\right\vert ^{2}=2\phi _{k}u_{k}g(\nabla u_{k},\nabla \phi _{k})+\frac{1%
}{2}\phi _{k}^{2}|\nabla u_{k}|^{2}+2u_{k}^{2}|\nabla \phi _{k}|^{2}
\end{equation*}%
to obtain 
\begin{equation*}
\int_{M}\Delta _{f}(u_{k})\phi _{k}^{2}u_{k}\left( e^{-f}d\mathrm{vol}%
_{g}\right) \leq -\frac{1}{2}\int_{B(x, r_k+1}\phi _{k}^{2}|\nabla
u_{k}|^{2}\left( e^{-f}d\mathrm{vol}_{g}\right) +2\int_{M}u_{k}^{2}|\nabla
\phi _{k}|^{2}\left( e^{-f}d\mathrm{vol}_{g}\right) .
\end{equation*}%
On the other hand 
\begin{eqnarray*}
\int_{M}\Delta _{f}(u_{k})\phi _{k}^{2}u_{k}\left( e^{-f}d\mathrm{vol}%
_{g}\right) &\geq &-\frac{1}{k}\int_{M}\phi _{k}^{2}u_{k}\left( e^{-f}d%
\mathrm{vol}_{g}\right) \\
&\geq &-\frac{2}{k}\int_{B(x,r_{k}+1)}u\left( e^{-f}d\mathrm{vol}_{g}\right)
\end{eqnarray*}%
So we have 
\begin{equation*}
\frac{1}{2}\int_{B(x,1)}|\nabla u_{k}|^{2}\left( e^{-f}d\mathrm{vol}%
_{g}\right) \leq \frac{8}{r_{k}^{2}}\int_{B(x,r_{k}+1)}u^{2}\left( e^{-f}d%
\mathrm{vol}_{g}\right) +\frac{2}{k}\int_{B(x,r_{k}+1)}u\left( e^{-f}d%
\mathrm{vol}_{g}\right) .
\end{equation*}

Note that, since the volume is finite, $u\in L^{2}(e^{-f}d\mathrm{vol}_{g})$
implies $u\in L^{1}(e^{-f}d\mathrm{vol}_{g})$ so the right hand side will go
to zero as $k\rightarrow \infty .$ Taking the limit and using that $u_{k}$
converge to $u$ in $W^{1,2}$ we obtain 
\begin{equation*}
\int_{B(x,1)}|\nabla u|^{2}(e^{-f}d\mathrm{vol}_{g})=0.
\end{equation*}%
Which implies $u$ is constant since it is continuous.
\end{proof}

\begin{remark}
One consequence of this theorem is that if a gradient shrinking soliton has $%
\mathrm{scal}\in L^2(e^{-f} d\mathrm{vol}_g)$ then either $\mathrm{scal}>0$
or the metric is flat (see \cite{Petersen-WylieI}).
\end{remark}

We can also apply the Yau-Naber Liouville theorem to obtain a strong minimum
principle for tensors. The strong minimum principle for tensors in the
parabolic setting were developed for the study of Ricci flow by Hamilton see 
\cite{HamiltonIII}.

\begin{theorem}[Tensor Minimum Principle]
\label{tensor_max} Let $\left( M,g,f\right) $ be a manifold with finite $f$%
-volume: $\int e^{-f}d\mathrm{vol}<\infty ,$ and $T$ a symmetric tensor on
some (tensor) bundle such that $\left\vert T\right\vert \in L^{2}(e^{-f}d%
\mathrm{vol}_{g})$ and 
\begin{equation*}
\Delta _{f}T=\lambda T+\Phi \left( T\right) ,\text{ where }g(\Phi
(T)(s),s)\leq 0\text{ and }\lambda >0,
\end{equation*}%
then $T$ is nonnegative and $\mathrm{ker}(T)$ is parallel.
\end{theorem}

Note that if $T$ and $\Phi \left( T\right) $ are nonnegative then the
Bochner formula shows that $T$ is parallel.

\begin{proof}
As long as $T$ is nonnegative and has constant rank Lemma \ref{T-lemma}
shows that $\mathrm{ker}(T)$ is parallel.

Denote the eigenvalues of $T$ by $\lambda _{1}\leq \lambda _{2}\leq \cdots .$
Let $s$ be a unit field such that $T\left( s\right) =\lambda _{1}s$ at $p$
otherwise extended by parallel translation along geodesics emanating from $p$%
. We can then calculate at $p\in M$%
\begin{eqnarray*}
\Delta _{f}\lambda _{1} &\leq &\Delta _{f}g\left( T\left( s\right) ,s\right)
\\
&=&g(\left( \Delta _{f}T\right) \left( s\right) ,s) \\
&=&\lambda g\left( T\left( s\right) ,s\right) +g\left( \Phi (T)(s),s\right)
\\
&\leq &\lambda \lambda _{1}
\end{eqnarray*}%
where the first inequality is in the barrier sense of Calabi (see \cite%
{Calabi}). Thus the first eigenvalue satisfies the differential inequality%
\begin{equation*}
\Delta _{f}\lambda _{1}\leq \lambda \lambda _{1}
\end{equation*}%
everywhere in the barrier sense. A similar analysis where we minimize over $%
k $ dimensional subspaces at a point shows that

\begin{equation*}
\Delta _{f}\left( \lambda _{1}+\cdots +\lambda _{k}\right) \leq \lambda
\left( \lambda _{1}+\cdots +\lambda _{k}\right)
\end{equation*}%
in the barrier sense.

To see that $T$ is nonnegative let $u=\min \left\{ \lambda _{1},0\right\} $,
then $u\leq 0$, 
\begin{equation*}
\Delta _{f}u\leq 0
\end{equation*}%
in the sense of barriers, and, since $\left\vert T\right\vert \in
L^{2}(e^{-f}d\mathrm{vol}_{g})$, so is $u$. The Yau-Naber Liouville Theorem
then implies $u$ is constant. In other words, either $\lambda _{1}\geq 0$ or 
$\lambda _{1}$ is constant and less than $0$. However, this last case is
impossible since if $\lambda _{1}$ is constant 
\begin{equation*}
0=\Delta _{f}\lambda _{1}\leq \lambda \lambda _{1}.
\end{equation*}

Knowing that $\lambda _{1}+\cdots +\lambda _{k}\geq 0,$ now allows us to
apply the strong minimum principle to show that, if $\lambda _{1}+\cdots
+\lambda _{k}$ vanishes at some point, then it vanishes everywhere (see \cite%
{McOwen} page 244). Since $\mathrm{dim}(\mathrm{ker}(T))$ is the largest $k$
such that $\lambda _{1}+\cdots +\lambda _{k}$ vanishes this shows that the
kernel is a distribution.
\end{proof}

We now apply the minimum principle to our formulas for the $f$-laplacian of
curvature. When $T=\mathcal{R}$ we have $\Phi \left( \mathcal{R}\right) =-2\left(%
\mathcal{R}^{2}+\mathcal{R}^{\#}\right)$. Since $\mathcal{R}^{2}$ is always
nonnegative we see from the minimum principle that the curvature operator of
a gradient shrinking soliton is nonnegative if and only if $\mathcal{R}^{\#}$
is nonnegative. In fact, by examining the proof of the minimum principle we
can also obtain the result alluded to in the introduction. As notation, let 
\begin{equation*}
\lambda _{1}\leq \lambda _{2}\leq \dots 
\end{equation*}%
be the ordering of the eigenvalues of the curvature operator.

\begin{corollary}
Let $(M,g,f)$ be a shrinking gradient Ricci soliton with $\lambda _{2}\geq 0$
and $\left\vert \mathcal{R}\right\vert \in L^{2}(e^{-f}d\mathrm{vol}_{g})$
then $\mathcal{R}\geq 0$, $\mathrm{ker}\mathcal{R}$ is parallel, and the
holonomy algebra 
\begin{equation*}
\mathfrak{hol}_{p}=\mathrm{im}\left( \mathcal{R}:\wedge
^{2}T_{p}M\rightarrow \wedge ^{2}T_{p}M\right) .
\end{equation*}
\end{corollary}

\begin{proof}
Fix a point $p$ and let $\phi _{1}$ be a parallel bi-vector such that 
\begin{equation*}
g\left( \mathcal{R}(\phi _{1}),\phi _{1}\right) =\lambda _{1}
\end{equation*}%
at $p$. Then, from the same argument as in the proof of the Tensor Minimum
principle, we obtain 
\begin{equation*}
\Delta _{f}\lambda _{1}\leq \lambda \lambda _{1}-g\left( \mathcal{R}%
^{\#}(\phi _{1}),\phi _{1}\right)
\end{equation*}

Let $\{\phi _{\alpha }\}$ be a basis of othonormal eigenvectors for $%
\mathcal{R}$. The structure constants of the Lie algebra are $C_{\alpha
\beta \gamma }=g\left( [\phi _{\alpha },\phi _{\beta }],\phi _{\gamma
}\right)$, which are fully anti-symmetric in $\alpha $, $\beta $, and $%
\gamma $. Then if $\lambda_2 \geq 0$, 
\begin{eqnarray*}
g\left( \mathcal{R}^{\#}(\phi _{1}),\phi _{1}\right) &=&\sum_{\alpha ,\beta
}\left( C_{1\alpha \beta }\right) ^{2}\lambda _{\alpha }\lambda _{\beta } \\
&=&\sum_{\alpha ,\beta \geq 2}\left( C_{1\alpha \beta }\right) ^{2}\lambda
_{\alpha }\lambda _{\beta } \\
&\geq &0.
\end{eqnarray*}%
Thus we see that $\lambda _{1}\geq 0.$ Next the tensor minimum principle can
be applied to see that $\mathrm{ker}\mathcal{R}$ is parallel. This shows in
turn that the orthogonal complement $\mathrm{im}\mathcal{R}$ is parallel.
The Ambrose-Singer theorem on holonomy then implies the last claim (see \cite%
{Besse}.)
\end{proof}

\begin{remark}
In dimension $3$ this implies that a gradient shrinking soliton with
nonnegative Ricci curvature has nonnegative curvature operator.
\end{remark}

\begin{remark}
There are simple examples of manifolds with $\lambda _{2}\geq 0$ that do not
admit $2$-nonnegative curvature operator metrics or even nonnegative Ricci
curvature metrics. Consider the product $N\times M$ where $N$ is a
negatively curved surface and $M$ a possibly one dimensional manifold with
nonnegative curvature operator. Then $\lambda _{1}<0$ and $\lambda _{2}=0$.
If $\mathrm{scal}_{M}>\left\vert \mathrm{scal}_{N}\right\vert ,$ then the
metric will also have positive scalar curvature.
\end{remark}

In the formula for $f$-Laplacian of the Ricci tensor we have $\Phi \left( 
\mathrm{Ric}\right) $ is $-2\mathcal{K}$ where 
\begin{equation*}
\mathcal{K}=\sum_{i=1}^{n}g(R(\cdot ,E_{i})(\mathrm{Ric}(E_{i})).
\end{equation*}%
If we let $\{E_{i}\}$ be a basis of eigenvectors for $\mathrm{Ric}$ with
eigenvalues $\rho _{i}$. Then 
\begin{equation*}
g(\mathcal{K}(Y),Y)=\sum_{i=1}^{n}\rho _{i}\mathrm{sec}(Y,E_{i})
\end{equation*}%
So that $\mathcal{K}\geq 0$ if $M$ has nonnegative (or nonpositive)
sectional curvature, or if $M$ is Einstein. The minimum principle gives the
following splitting theorem for shrinking solitons with $\mathcal{K}\geq 0$.
This is the soliton version of a result of B\"{o}hm and Wilking \cite%
{Bohm-Wilking2007} which states that any compact manifold with nonnegative
sectional curvature and finite fundamental group flows in a short time under
the Ricci flow to a metric with positive Ricci curvature.

\begin{corollary}
\label{splittingII}Let $(M,g,f)$ be a shrinking gradient Ricci soliton with $%
\mathcal{K}\geq 0$ and $\left\vert \mathrm{Ric}\right\vert \in L^{2}(e^{-f}d%
\mathrm{vol}_{g})$ then $\tilde{M}=N\times \mathbb{R}^{k}$ where $N$ has
positive Ricci curvature. In particular, a compact shrinking soliton with $%
\mathcal{K}\geq 0$ has positive Ricci curvature.
\end{corollary}

\begin{proof}
The minimum principle and the de Rham splitting theorem show that $\tilde{M}%
=N\times F,$ where $N$ has positive Ricci curvature and $F$ is Ricci flat.
From \cite{Petersen-WylieII} we get that both $N\ $and $F$ are gradient
solitons. Finally Ricci flat solitons are Gaussians, thus proving the
corollary. The last bit about compact manifolds follows from the fact that
shrinking solitons have finite volume and hence finite fundamental group.
\end{proof}

We now turn our attention to the proof of Theorem \ref{shrinker-ndim}. We
assume that $M$ is a non-flat, gradient shrinking soliton with $W=0$ and $%
\left\vert \mathrm{Ric}\right\vert \in L^{2}(e^{-f}d\mathrm{vol}_{g}).$
Recall that when $W=0$ 
\begin{eqnarray*}
\Delta _{f}\mathrm{Ric} &=&2\lambda \mathrm{Ric}-\frac{2n\mathrm{scal}}{%
\left( n-1\right) \left( n-2\right) }\mathrm{Ric}+\frac{4}{n-2}\mathrm{Ric}%
^{2} \\
&&-\frac{2}{\left( n-2\right) }\left( \left\vert \mathrm{Ric}\right\vert
^{2}-\frac{\mathrm{scal}^{2}}{n-1}\right) I
\end{eqnarray*}%
Let $\rho _{1}\leq \cdots \leq \rho _{n}$ be the eigenvalues of $\mathrm{Ric}
$ and $E$ a unit field such that $\mathrm{Ric}\left( E\right) =\rho _{1}E$
at $p\in M$ and extend it to be parallel along geodesics emmanating from $p.$
Clearly $\rho _{1}\leq \mathrm{Ric}\left( E,E\right) $ with equality at $p.$
Calculating at $p$ we have 
\begin{eqnarray*}
\Delta _{f}\left( \rho _{1}\right) &\leq &\Delta _{f}\mathrm{Ric}\left(
E,E\right) \\
&=&\left( \Delta _{f}\mathrm{Ric}\right) \left( E,E\right) \\
&=&2\lambda \rho _{1}-\frac{2n\mathrm{scal}}{\left( n-1\right) \left(
n-2\right) }\rho _{1}+\frac{4}{n-2}\rho _{1}^{2}-\frac{2}{n-2}\left(
\left\vert \mathrm{Ric}\right\vert ^{2}-\frac{\mathrm{scal}^{2}}{n-1}\right)
\end{eqnarray*}%
where $\Delta _{f}\left( \rho _{1}\right) $ is interpreted as being in the
upper barrier sense of \cite{Calabi}.

The ratio $\frac{\rho _{1}}{\mathrm{scal}}$ then satisfies%
\begin{eqnarray*}
\Delta _{h}\left( \frac{\rho _{1}}{\mathrm{scal}}\right) &\leq &2\phi , \\
h &=&f-\log \left( \mathrm{scal}^{2}\right) ,
\end{eqnarray*}%
where 
\begin{eqnarray*}
\phi &=&\frac{\rho _{1}^{2}\left( n\rho _{1}-\mathrm{scal}\right) }{\left(
n-1\right) \mathrm{scal}^{2}} \\
&&+\frac{\left( (n-2)\rho _{1}-\mathrm{scal}\right) \left(
(n-1)\sum_{j=2}^{n}(\rho _{j})^{2}-\left( \sum_{j=2}^{n}\rho _{j}\right)
^{2}\right) }{(n-1)(n-2)\mathrm{scal}^{2}}
\end{eqnarray*}%
which is clearly nonpositive.

We now have $\frac{\rho_1}{\mathrm{scal}} \leq 1$ and $\Delta _{h}\left( 
\frac{\rho _{1}}{\mathrm{scal}}\right) \leq 0$, so to apply the Yau-Naber
Liouville Theorem, we must show the measure is finite and the function is in 
$L^2(e^{-h} d \mathrm{vol}_g$). This is clear from Lemma \ref{vol-comp}
because 
\begin{equation*}
\int_M e^{-h} d\mathrm{vol}_g = \int_M \mathrm{scal}^2 e^{-f}d\mathrm{vol}_g
<\infty
\end{equation*}
and 
\begin{equation*}
\int_M \left(\frac{\rho_1}{\mathrm{scal}}\right)^2 e^{-h}d\mathrm{vol}_g =
\int_M \left(\rho_1\right)^2 e^{-f} d\mathrm{vol}_g < \infty
\end{equation*}
Thus $\frac{\rho _{1}}{\mathrm{scal}}$is constant. In particular, $\phi $
must vanish and

\begin{eqnarray*}
\rho _{1}^{2}\left( n\rho _{1}-\mathrm{scal}\right) &=&0 \\
\left( (n-2)\rho _{1}-\mathrm{scal}\right) \left( (n-1)\sum_{j=2}^{n}(\rho
_{j})^{2}-\left( \sum_{j=2}^{n}\rho _{j}\right) ^{2}\right) &=&0
\end{eqnarray*}

The first equation tells us that either $\rho _{1}=0$ or $M$ is Einstein.
When $\rho _{1}=0$ the second equation and the Cauchy-Schwarz inequality
tells us that 
\begin{equation*}
\rho _{2}=\rho _{3}=\cdots =\rho _{n}=\frac{\mathrm{scal}}{n-1}>0.
\end{equation*}%
Then by Corollary \ref{splittingI} the universal cover of $M$ splits $\tilde{%
M}=N\times \mathbb{R}$ where $N$ is again a shrinking gradient soliton with
a Ricci tensor that has only one eigenvalue. When $n=3$ Hamilton's
classification of surface solitons (see Appendix) then shows that $N$ is the
standard sphere, while if $n>3$ Schur's lemma shows that $N$ is Einstein.

\appendix

\section{Surface gradient solitons}

In the literature the classification of shrinking surface solitons is
usually stated for metrics with bounded curvature. However, we have used
this classification under the weaker condition $\mathrm{scal}\in
L^{2}(e^{-f}d\mathrm{vol}_{g})$. In this appendix we verify that the
classification still holds in this case.

We consider warped product metrics, a slightly larger class of metrics than
rotationally symmetric ones.

\begin{definition}
A Riemannian metric $(M,g)$ on either $\mathbb{R}^{n}$, $S^{n}$, or $N\times 
\mathbb{R}$ is a warped product if it can be written as 
\begin{equation*}
g=dr^{2}+h^{2}(r)g_{0}.
\end{equation*}%
When $M=\mathbb{R}^{n}$, we assume that $h(0)=0$ and $g_{0}$ is the standard
metric on the sphere. When $M=S^{n}$, we require $h(0)=h(r_{0})=0$ and $%
g_{0} $ is the standard metric on the sphere.
\end{definition}

There is a very simple Obata-type characterization of warped product metrics
found in \cite{Cheeger-Colding}.

\begin{theorem}[Cheeger-Colding]
A Riemannian manifold $(M,g)$ is a warped product if and only if there is a
nontrivial function $f$ such that 
\begin{equation*}
\mathrm{Hess}f=\mu g
\end{equation*}%
for some function $\mu :M\rightarrow \mathbb{R}$.
\end{theorem}

\begin{proof}
If $g=dr^{2}+h^{2}(r)g_{0}$ simply let $f=\int h\left( r\right) dr.$

Conversely, we see that $f$ is rectifiable (see \cite{Petersen-WylieII}) as%
\begin{equation*}
D_{X}\frac{1}{2}\left\vert \nabla f\right\vert ^{2}=\mathrm{Hess}f\left(
X,\nabla f\right) =\mu g\left( X,\nabla f\right)
\end{equation*}%
Showing that $\left\vert \nabla f\right\vert $ is constant on level sets of $%
f.$ Let $N$ be a nondegenerate level set of $f$, $g_{0}$ the metric
restricted to this level set, and $r$ the signed distance to $N$ defined to
that $\nabla r$ and $\nabla f$ point in the same direction. Then $f=f\left(
r\right) $%
\begin{eqnarray*}
\nabla f &=&f^{\prime }\nabla r, \\
\mathrm{Hess}f &=&f^{\prime \prime }dr^{2}+f^{\prime }\mathrm{Hess}r.
\end{eqnarray*}%
This shows that $\mu =f^{\prime \prime }$ and that%
\begin{equation*}
\mathrm{Hess}r=\frac{f^{\prime \prime }}{f^{\prime }}g
\end{equation*}%
on the orthogonal complement of $\nabla r.$ Thus $g=dr^{2}+\left( cf^{\prime
}\right) ^{2}g_{0}$ where $cf^{\prime }\left( 0\right) =1.$
\end{proof}

\begin{remark}
This theorem indicates that any flow that preserves conformal classes has
the property that the corresponding gradient solitons must be rotationally
symmetric. This then makes it possible to classify all complete gradient
solitons for such flows.
\end{remark}

\begin{corollary}
Any surface gradient Ricci soliton is a warped product.
\end{corollary}

\begin{proof}
Simply use that 
\begin{equation*}
\mathrm{Ric}=\frac{\mathrm{scal}}{2}g.
\end{equation*}
\end{proof}

So the problem of finding surface gradient solitons is reduced to
determining which functions $h(r)$ give a soliton. For example, Hamilton's
cigar is obtained by taking $h(r) = \tanh(r)$ and is the unique (up to
scaling) non-compact steady gradient soliton surface with positive
curvature. For a non-trivial example of an expanding surface gradient
soliton see (\cite{Chow-Lu-Ni}, p. 164-167).

Now suppose we have a non-flat shrinking soliton on a surface with $\mathrm{%
scal}\in L^2(e^{-f}d \mathrm{vol}_g)$. As we have seen this implies $\mathrm{%
scal}>0$. Moreover, since we are on a surface, the Ricci curvature is
positive so Proposition 1.1 in \cite{Ni} implies the scalar curvature is
bounded away from zero. Thus $M$ is compact. Then, since it is also a warped
product, $M$ must be a rotationally symmetric metric on the sphere. Chen,
Lu, and Tian show that this implies $M$ is a round sphere \cite{Chen-Lu-Tian}%
.

\end{document}